\documentclass[a4paper,12pt,reqno]{amsart}

\usepackage{amsmath,amsfonts,amsthm,amssymb}
\usepackage{color}

\numberwithin{equation}{section}
\newtheorem{theorem}{Theorem}[section]
\newtheorem{lemma}[theorem]{Lemma}

\theoremstyle{definition}
\newtheorem{definition}[theorem]{Definition}

\newtheorem{remark}[theorem]{Remark}

\newcommand{\dist}{\operatorname{dist}}

\newcommand{\Ca}{\operatorname{cap}}
\newcommand{\con}{\operatorname{con}}

\def\kint_#1{\mathchoice%
          {\mathop{\kern 0.2em\vrule width 0.6em height 0.69678ex depth -0.58065ex
                  \kern -0.8em \intop}\nolimits_{\kern -0.4em#1}}%
          {\mathop{\kern 0.1em\vrule width 0.5em height 0.69678ex depth -0.60387ex
                  \kern -0.6em \intop}\nolimits_{#1}}%
          {\mathop{\kern 0.1em\vrule width 0.5em height 0.69678ex depth -0.60387ex
                  \kern -0.6em \intop}\nolimits_{#1}}%
          {\mathop{\kern 0.1em\vrule width 0.5em height 0.69678ex depth -0.60387ex
                  \kern -0.6em \intop}\nolimits_{#1}}}
\def\vintslides_#1{\mathchoice%
          {\mathop{\kern 0.1em\vrule width 0.5em height 0.697ex depth -0.581ex
                  \kern -0.6em \intop}\nolimits_{\kern -0.4em#1}}%
          {\mathop{\kern 0.1em\vrule width 0.3em height 0.697ex depth -0.604ex
                  \kern -0.4em \intop}\nolimits_{#1}}%
          {\mathop{\kern 0.1em\vrule width 0.3em height 0.697ex depth -0.604ex
                  \kern -0.4em \intop}\nolimits_{#1}}%
          {\mathop{\kern 0.1em\vrule width 0.3em height 0.697ex depth -0.604ex
                  \kern -0.4em \intop}\nolimits_{#1}}}

\def\mvint_#1{\mathchoice
          {\mathop{\vrule width 6pt height 3 pt depth -2.5pt
                  \kern -8pt \intop}\nolimits_{\hspace{-1ex}#1}}%
          {\mathop{\vrule width 5pt height 3 pt depth -2.6pt
                  \kern -6pt \intop}\nolimits_{\hspace{-1ex}#1}}%
          {\mathop{\vrule width 5pt height 3 pt depth -2.6pt
                  \kern -6pt \intop}\nolimits_{\hspace{-1ex}#1}}%
          {\mathop{\vrule width 5pt height 3 pt depth -2.6pt
                  \kern -6pt \intop}\nolimits_{\hspace{-1ex}#1}}}
\newcommand{\dd}{\,\textrm{d}}

\newcommand{\R}{\mathbb{R}}
\newcommand{\Z}{\mathbb{Z}}

\begin{document}

\author{Juha Kinnunen and Riikka Korte}
\title[Characterizations of Sobolev inequalities]{Characterizations of Sobolev inequalities on metric spaces}

\begin{abstract}
We present isocapacitary characterizations of Sobolev 
inequalities in very general metric measure spaces. 
\end{abstract}
\maketitle

\section{Introduction}
There is a well known connection between the isoperimetric and
Sobolev inequalities. By the isoperimetric inequality, we have
\begin{equation}\label{isoperineq}
|E|^{(n-1)/n}\le c(n)\mathcal H^{n-1}(\partial E),
\end{equation}
where $E$ is a smooth enough subset of $\R^n$, 
$|E|$ is the Lebesgue measure and $\mathcal H^{n-1}$ is the 
$(n-1)$-dimensional Hausdorff measure. The constant $c(n)$
is chosen so that \eqref{isoperineq} becomes an equality
when $E$ is a ball. 
The Sobolev inequality states that
\begin{equation}\label{sobolev}
\bigg(\int_{\R^n}|u|^{n/(n-1)}\dd x\bigg)^{(n-1)/n}
\le c(n)\int_{\R^n}|\nabla u|\dd x
\end{equation}
for every $u\in C^\infty_0(\R^n)$.  The smallest constant in
\eqref{sobolev} is the same as the constant in \eqref{isoperineq}. The
Sobolev inequality follows from the isoperimetric inequality through
the co-area formula. On the other hand, the isoperimetric inequality can
be deduced from the Sobolev inequality, see for
example~\cite{evansgariepy}. This shows that the isoperimetric and
Sobolev inequalities are different aspects of the same phenomenon.
Originally this observation is due to Federer and Fleming
\cite{fedfle} and Maz$'$ya \cite{mazja60}.

When the gradient is integrable to a power which is greater than one,
the isoperimetric inequality has to be replaced with an isocapacitary
inequality. When the exponent is one, capacity and Hausdorff content
are equivalent and hence it does not matter which one we choose. In
this case due to the boxing inequality, it is enough to have the
isocapacitary inequality for balls instead of all sets.  This elegant
approach to Sobolev inequalities is due to Maz$'$ya, see \cite{mazja85}
and~\cite{mazja03}. Usually this characterization leads to
descriptions of the best possible constants in Sobolev inequalities.
However, the aim of the present work is not so much to study best
possible constants but rather study necessary and sufficient
conditions for Sobolev--Poincar\' e inequalities in a general metric
space context. In weighted Euclidean spaces, these characterizations
have been studied in~\cite{Turesson}.

Rather standard assumptions in analysis on metric measure spaces
include a doubling condition for the measure and validity of some kind
of Sobolev--Poincar\' e inequality. Despite the fact that plenty of
analysis has been done in this general context, very little is known
about the basic assumptions.  Several necessary conditions are known,
but unfortunately only few sufficient conditions are available so far.
On Riemannian manifolds, Grigor'yan and Saloff-Coste observed that the
doubling condition and the Poincar\' e inequality are not only
sufficient but also necessary conditions for a scale invariant
parabolic Harnack principle for the heat equation, see \cite{SaCo1},
\cite{SaCo2} and \cite{Gri}.  It also known that Maz$'$ya type
characterizations of Sobolev inequalities are available on Riemannian
manifolds.  The purpose of this work is to show
that this is also the case on a very general metric measure spaces.

\section{Preliminaries}
We assume that $X=(X,d,\mu)$ is a metric measure space equipped with a 
metric $d$ and a Borel regular outer measure $\mu$ such that $0<\mu(B)<\infty$ 
for all balls $B=B(x,r)=\{y\in X\,:\, d(x,y)<r\}$. In what follows, 
$\Omega$ stands for an open bounded subset of $X$ unless otherwise stated. 
The measure $\mu$ is said to be \emph{doubling} if there exists a constant 
$c_D\geq 1$, called the \emph{doubling constant}, such that
\[
\mu(B(x,2r))\leq c_D\mu(B(x,r))
\]
for all $x\in X$ and $r>0$.

In this paper, a \emph{path} in $X$ 
is a rectifiable non-constant continuous mapping from a
compact interval to $X$. 
A path can thus be parameterized by arc length. 

By saying that a condition holds for \emph{$p$-almost every path}
with $1\le p<\infty$,
we mean that it fails only for a path family with zero $p$-modulus. 
A family $\Gamma$ of curves is of zero $p$-modulus if there is a 
non-negative Borel measurable function $\rho\in L^p(X)$ such that 
for all curves $\gamma\in\Gamma$, the path integral $\int_\gamma \rho\dd s$ is infinite.

A nonnegative Borel function $g$ on $X$ is an \emph{upper gradient} 
of an extended real valued function $u$
on $X$ if for all paths $\gamma$ 
joining points $x$ and $y$ in $X$ we have
\begin{equation} \label{ug-cond}
|u(x)-u(y)|\le \int_\gamma g\dd s,
\end{equation}
whenever both $u(x)$ and $u(y)$ are finite, and 
$\int_\gamma g\dd s=\infty $ otherwise.
If $g$ is a nonnegative measurable function on $X$
and if (\ref{ug-cond}) holds for $p$-almost every path,
then $g$ is a \emph{$p$-weak upper gradient} of~$u$.

Let $1\le p<\infty$.
If $u$ is a function that is integrable to power $p$ in $X$, let
$$
\|u\|_{N^{1,p}(X)} 
= \Big(\int_X |u|^p\dd\mu 
+\inf_g\int_X g^p \dd\mu\Big)^{1/p},
$$
where the infimum is taken over all $p$-weak upper gradients of $u$.
The \emph{Newtonian space} on $X$ is the quotient space
$$
N^{1,p}(X)=\{u:\|u\|_{N^{1,p}(X)}<\infty\}/{\sim},
$$
where  $u \sim v$ if and only if $\|u-v\|_{N^{1,p}(X)}=0$.
For properties of Newtonian spaces, we refer to \cite{Sh-rev}.

Let $E$ be a subset of $\Omega$. We write $u\in\mathcal A(E,\Omega)$ 
if $u_{|E}=1$ and $u_{|X\setminus \Omega}=0$.

\begin{definition} Let $E\subset \Omega$. 
  The \emph{$p$--capacity} of $E$ with respect to $\Omega$ is
  \[
  \Ca_p(E,\Omega)=\inf\int_\Omega g_u^p\dd\mu,
  \]
  where the infimum is taken over all continuous functions
  $u\in\mathcal A(E,\Omega)$ with $p$--weak upper gradients $g_u$. If
  there are no such functions, then $\Ca_p(E,\Omega)=\infty$.
\end{definition}

We say that a property regarding points in $X$ 
holds \emph{$p$--quasieverywhere} ($p$--q.e.)
if the set of points  for which the property does not hold 
has capacity zero. 

To be able to compare the boundary values of Newtonian functions, we
need a Newtonian space with zero boundary values.  Let $E$ be a
measurable subset of $X$. The \emph{Newtonian space with zero boundary
  values} is the space
\[
N^{1,p}_0(E)=\{u|_{E} : u \in N^{1,p}(X) \text{ and }
        u=0 \text{ $p$--q.e. on } X \setminus E\}.
\]
Note that if $\Ca_p(X \setminus E,X) = 0$, then $N^{1,p}_0(E) =
N^{1,p}(X)$.  The space $N^{1,p}_0(E)$ equipped with the norm
inherited from $N^{1,p}(X)$ is a Banach space, see Theorem~4.4
in~\cite{Sh-harm}.

\begin{definition}
  We say that $X$ supports a \emph{weak $(1,p)$--Poincar\'e
    inequality} if there exist constants $c_p>0$ and $\tau\geq 1$ such
  that for all balls $B(x,r)$ of $X$, all integrable functions
  $u$ on $X$ and all upper gradients $g_u$ of $u$, we have
\begin{equation}\label{eqn:poincare}
  \kint_{B(x,r)}|u-u_{B(x,r)}|\dd\mu\leq c_pr
  \left(
    \kint_{B(x,\tau r)}g_u^p \dd\mu
  \right)^{1/p},
\end{equation}
where
\[
u_B=\kint_B u\dd\mu=\frac{1}{\mu(B)}\int_Bu\dd\mu.
\]
\end{definition}

\section{Functions with zero boundary values}

In this section, we give necessary and sufficient conditions for 
Sobolev inequalities of type
\[
\left(\int_\Omega |u|^q\dd\nu\right)^{1/q}
\leq c_S\left(\int_\Omega g_u^p\dd\mu\right)^{1/p},
\]
where the constant $c_S$ is independent of 
$u\in N^{1,p}_0(\Omega)\cap C(\Omega)$, 
and $\mu$ and $\nu$ are Borel regular outer measures.
We consider two ranges of indices separately.

We have chosen to study continuous Newtonian functions, 
but our arguments do not depend on the choice of the function space.
For example, similar results also hold for all Newtonian functions and 
Lipschitz functions if the definition of $p$--capacity is adjusted accordingly.

\begin{remark}
In the standard versions of Poincar\'e inequality, the inequality
depends on the diameter of the set. Therefore the constant
$c_S$ also depends strongly on the diameter of the set in many cases.
\end{remark}

\subsection{ The case $1\leq p\leq q<\infty$}
Let $u\,:\,\Omega\rightarrow [-\infty,\infty]$ be a $\mu$--measurable function.
By the well-known Cavalieri principle
\[
\int_\Omega|u|^p\dd\mu
 =p\int_0^\infty \lambda^{p-1}\mu(E_\lambda)\dd\lambda,
\]
where
\[
E_\lambda=\{x\in\Omega\,:\, |u(x)|>\lambda\}.
\]
The following simple integral inequality will
be useful for us. 
Notice that the equality occurs when $p=q$.

\begin{lemma}\label{lemma:qpey}
  If $u\,:\,\Omega\rightarrow [-\infty,\infty]$ is $\mu$--measurable 
  and $0<p\leq q<\infty$, then
  \begin{equation}\label{eqn:qpey}
    \left(
      \int_\Omega |u|^q\dd\mu
    \right)^{1/q}
    \leq\left(
      p\int_0^\infty \lambda^{p-1}\mu(E_\lambda)^{p/q}\dd\lambda
    \right)^{1/p}.
  \end{equation}
\end{lemma}

\begin{proof}  
  We have
  \[
  \left(
    \int_X
    \left(
      |u|^p
    \right)^{q/p}\dd\mu
  \right)^{p/q}
  =\sup_{\|f\|_s\leq 1}\int_X|u|^pf\dd\mu,
  \]
  where $s=q/(q-p)$ is the H\"older conjugate of $q/p$. 
  By $\|f\|_s$ we denote the  $L^s(\mu)$--norm of $f$.
  Define a measure $\widetilde\mu$ as 
  \[
  \widetilde\mu(A)=\int_A|f|\dd\mu
  \]
  for every $\mu$--measurable set $A\subset X$.
  If $\mu(E)>0$ and $\|f\|_s\le1$, we conclude that
  \[
\begin{split}
  \widetilde\mu(E)&=\int_E|f|\dd\mu\leq\mu(E)^{1-1/s}
  \left(
    \int_E|f|^s\dd\mu
  \right)^{1/s}
\\
  &\leq \mu(E)^{1-1/s}=\mu(E)^{p/q}.
\end{split} 
 \]
  Here we used the H\"older inequality. Hence
  \begin{equation*}
    \begin{split}
     \int_X|u|^pf\dd\mu &\leq\int_X|u|^p\dd\widetilde\mu 
      =p \int_0^\infty \lambda^{p-1}\widetilde\mu(E_\lambda)\dd\lambda
      \\
&\leq p\int_0^\infty\lambda^{p-1}\mu(E_\lambda)^{p/q} \dd\lambda.
    \end{split}
  \end{equation*}
  Taking supremum over all functions $f$ with $\|f\|_s\leq 1$ completes the proof.
\end{proof}

Next we prove a strong type inequality for the capacity.  When $p=1$
the obtained estimate reduces to the co-area formula.  The proof is
based on a general truncation argument, see page 110 in
\cite{mazja85}. A similar argument has been used for example
in~\cite{BCLS},~\cite{SaCo1} and~\cite{HK}.

\begin{lemma}\label{lemma:capaint}
  Let $u\in N^{1,p}_0(\Omega)\cap C(\Omega)$ and $1\leq p<\infty$. Then 
  \[
  \int_0^\infty \lambda^{p-1}\Ca_p(E_\lambda,\Omega)\dd\lambda 
  \leq 2^{2p-1}\int_\Omega g_u^p\dd\mu,
  \]
  where $g_u$ is a $p$-weak upper gradient of $u$.
\end{lemma}
\begin{proof}
A straightforward calculation shows that
  \[
    \begin{split}
      &\int_0^\infty \lambda^{p-1}\Ca_p(E_\lambda,\Omega)\dd\lambda\\
      &=\sum_{j=-\infty}^\infty
      \int_{2^{j-1}}^{2^j}\lambda^{p-1}\Ca_p(E_\lambda,\Omega)\dd\lambda\\
      &\leq\sum_{j=-\infty}^\infty 
      (2^j-2^{j-1})2^{j(p-1)}\Ca_p(E_{2^{j-1}},\Omega)\\
      &=\frac{1}{2}\sum_{j=-\infty}^\infty 2^{jp}\Ca_p(E_{2^{j-1}},\Omega)\\
      &=2^{p-1}\sum_{j=-\infty}^\infty 2^{jp}\Ca_p(E_{2^j},\Omega).
    \end{split}
  \] 
  Let
  \begin{equation*}
    u_j=
    \begin{cases}
      1, & \text{if}\quad u\geq 2^j,\\
      2^{1-j}|u|-1, & \text{if}\quad 2^{j-1}<u<2^j,\\
      0,&\text{if}\quad u\leq 2^{j-1}. 
    \end{cases}
  \end{equation*}
  Then $u_j\in \mathcal A(E_{2^j},\Omega)$. This implies that
  \[
  \Ca_p(E_{2^j},\Omega)
  \leq 2^{p(1-j)}\int_{E_{2^{j-1}}\setminus E_{2^{j}}} g_u^p\dd\mu
  \]
  and consequently
  \begin{equation*}
    \begin{split}
      \sum_{j=-\infty}^\infty 2^{jp}\Ca_p(E_{2^j},\Omega)
      &\leq \sum_{j=-\infty}^\infty 
      2^{jp+ p(1-j)}\int_{E_{2^{j-1}}\setminus E_{2^{j}}} g_u^p\dd\mu\\
      &\leq 2^p\int_\Omega g_u^p\dd\mu.
    \end{split}
  \end{equation*}
  The claim follows from this.
\end{proof}

The following result gives a necessary and sufficient condition
for a Sobolev inequality in terms of an isocapacitary inequality. 
This is a metric space version of a corollary on page 113
of \cite{mazja85}.  

\begin{remark} We do not need the doubling condition in this theorem.
\end{remark}

\begin{theorem}\label{theorem:p_pienempi}
  Suppose that $1\leq p\leq q<\infty$.
  \begin{enumerate}
  \item[(i)] If there is a constant $\gamma$ such that
    \begin{equation}\label{eqn:isoperi}
      \nu(E)^{p/q}\leq\gamma\Ca_p(E,\Omega)
    \end{equation}
    for every $E\subset\Omega$, then
    \begin{equation}\label{eqn:sobolev}
      \left(
        \int_\Omega |u|^q\dd\nu
      \right)^{1/q}
      \leq c_S\left(\int_\Omega g_u^p\dd\mu\right)^{1/p}
    \end{equation}
    for every $u\in N^{1,p}_0(\Omega)\cap C(\Omega)$ with $c_S$ depending only on $\gamma$ and $p$.
  \item[(ii)] If~\eqref{eqn:sobolev} holds for every $u\in N^{1,p}_0(\Omega)\cap C(\Omega)$ 
    and if the constant $c_S$ is independent of $u$, then~\eqref{eqn:isoperi} 
    holds for every $E\subset\Omega$ with $\gamma= c_S$.
  \end{enumerate}
\end{theorem}

\begin{proof}
  (i)
    By Lemma~\ref{lemma:qpey}, \eqref{eqn:isoperi} and 
    Lemma~\ref{lemma:capaint}, we obtain
    \begin{equation*}
      \begin{split}
        \left(\int_\Omega |u|^q\dd\nu\right)^{1/q}&
        \leq \left(p\int_0^\infty \lambda^{p-1}\nu(E_\lambda)^{p/q}
          \dd\lambda\right)^{1/p}\\
        &\leq\left(\gamma p\int_0^\infty\lambda^{p-1}
          \Ca_p(E_\lambda,\Omega)\dd\lambda\right)^{1/p}\\
        &\leq\left(\gamma p2^{2p-1}\right)^{1/p}
        \left(\int_\Omega g_u^p\dd\mu\right)^{1/p}.
      \end{split}
    \end{equation*}
(ii) 
    If $u\in\mathcal A(E,\Omega)$ is continuous, then by~\eqref{eqn:sobolev}, we have
    \begin{equation*}
      \nu(E)^{1/q}
      \leq \left(\int_\Omega |u|^q\dd\nu\right)^{1/q}
      \leq c_S\left(\int_\Omega g_u^p\dd\mu\right)^{1/p}.
    \end{equation*}
    The claim follows by taking the infimum on the right-hand side.
\end{proof}

\begin{remark}
The previous theorem gives a necessary and sufficient condition
for the Hardy inequality
\begin{equation}\label{hardy}
\int_\Omega\left(\frac{|u(x)|}{\dist(x,X\setminus\Omega)}\right)^p\dd\mu
\le c_H\int_\Omega g_u^p\dd\mu,
\end{equation}
where the constant $c_H$ is independent of $u\in N^{1,p}_0(\Omega)\cap
C(\Omega)$.  Indeed, \eqref{hardy} holds if and only if
\[
\int_E\frac1{\dist(x,X\setminus\Omega)^p}\dd\mu
\leq\gamma\Ca_p(E,\Omega)
\]
for every $E\subset\Omega$.  Thus Theorem~\ref{theorem:p_pienempi} is a
generalization of Theorem 4.1 in~\cite{KS}.  
In the metric space context, the Hardy
inequality has also been studied in \cite{BMS}.
\end{remark}

\subsection{ The case $p=1$}

When $p=1$ and $\Omega=X$, the isocapacitary inequalities reduce to
isoperimetric inequalities.  Moreover, in this case we can improve
Theorem~\ref{theorem:p_pienempi} under the additional assumptions that
the measure is doubling and the space supports a Poincar\'e
inequality.  Indeed, it is enough that condition~\eqref{eqn:isoperi}
is satisfied for all balls. To prove that, we will need equivalence of
the capacity of order one and the Hausdorff content of co-dimension
one
\[
\mathcal H^h_\infty(K)=
\inf\left\{\sum_{i=1}^\infty\frac{\mu(B(x_i,r_i))}{r_i}\,:\, K\subset\bigcup_{i=1}^\infty B(x_i,r_i)\right\}.
\]
\begin{theorem}\label{theorem:kontentti_kapasiteetti} Let $X$ be a
  complete metric space with a doubling measure $\mu$. Suppose that
  $X$ supports a weak $(1,1)$--Poincar\'e inequality.  Let $K$ be a
  compact subset of $X$. Then
\[
\frac{1}{c}\Ca_1(K)\leq \mathcal H^h_\infty(K)\leq c \Ca_1(K),
\]
where $c$ depends only on the doubling constant and the constants in
the weak $(1,1)$--Poincar\'e inequality.
\end{theorem}

The proof is based on co-area formula and a metric space version of
so--called boxing inequality.  For more details, see~\cite{kkst06}. A
similar result has been studied in M\"ak\"al\"ainen~\cite{makalainen}.

\begin{theorem}\label{theorem:pallot}
  Let $X$ be a complete metric space with a doubling measure $\mu$. Suppose
  that $X$ supports a weak $(1,1)$--Poincar\'e inequality. Suppose
  that $1\leq q<\infty$.  If there is a constant $\gamma$ such that
  \begin{equation}\label{eqn:isoperi_pallo}
    \nu(B)^{1/q}\leq\gamma\Ca_1(B,X)
  \end{equation}
  for every ball $B\subset X$, then
  \begin{equation}\label{eqn:sobolev_pallo}
    \left(
      \int_X |u|^q\dd\nu
    \right)^{1/q}
    \leq c_S\int_X g_u\dd\mu,
  \end{equation}
  where $c_S$ is independent of $u\in N^{1,1}(X)\cap C_0(X)$.
\end{theorem}

\begin{proof}
  First we prove that if the space satisfies~\eqref{eqn:isoperi_pallo}
  for all balls in $X$, then it satisfies the same condition for all
  compact sets with a different constant.

  Let $K\subset X$ be compact, $\varepsilon>0$ and
  $\{B(x_i,r_i)\}_{i=1}^\infty$ be a covering of $K$ such that
  \[
  \mathcal H^h_\infty(K)\geq
  \sum_{i=1}^\infty\frac{\mu(B(x_i,r_i))}{r_i}-\varepsilon.
  \]
  Since $q\geq1$, we have
  \[
  \nu(K)^{1/q}\leq \sum_{i=1}^\infty\nu(B(x_i,r_i))^{1/q}.
  \]
  Because 
  \[
  u_i(x)=\left(1-\dist(x,B(x_i,r_i))/r_i\right)_+
  \] 
  belongs to $ \mathcal A(B(x_i,r_i),X)$, and
  $g_i=\chi_{B(x_i,2r_i)}/r_i$ is an upper gradient of $u_i$, we have
  \[
  \Ca_1(B(x_i,r_i))\leq \int_X g_i\dd\mu\leq
  c_D\frac{\mu(B(x_i,r_i))}{r_i}.
  \]
  By combining the above estimates and~\eqref{eqn:isoperi_pallo}, we
  conclude
  \begin{equation*}
    \begin{split}
      \nu(K)^{1/q} 
      &\leq \sum_{i=1}^\infty \nu(B(x_i,r_i))^{1/q}\\
      &\leq \gamma\sum_{i=1}^\infty\Ca_1(B(x_i,r_i))\\
      &\leq \gamma c_D\sum_{i=1}^\infty \frac{\mu(B(x_i,r_i))}{r_i}\\
      &\leq \gamma c_D (\mathcal H^h_\infty(K)+\varepsilon).
    \end{split}
  \end{equation*}
  The claim follows by
  Theorem~\ref{theorem:kontentti_kapasiteetti} as
  $\varepsilon\rightarrow 0$.

  Now the theorem follows as in the proof of
  Theorem~\ref{theorem:p_pienempi}. Note that since $u$ has compact
  support and is continuous, we can as well consider compact level
  sets $\{|u|\geq t\}$ instead of open sets.
\end{proof}

\subsection{ The case $1\leq q<p<\infty$}
In the case $1\leq q<p<\infty$, the isocapacitary inequality
takes a different form. Let $E_j$, $j=-N,-N+1,\ldots, N,N+1$ 
be such that $ E_j\subset \Omega$ and $ E_j\subset E_{j+1}$ for $j=-N,-N+1,\ldots, N$. 
We define
\begin{equation}\label{eqn:gamma}
  \gamma=\sup\left[
    \sum_{j=-N}^{N}
    \left(\frac{\nu(E_j)^{p/q}}{\Ca_p(E_j,E_{j+1})}
    \right)^{q/(p-q)}
  \right]^{(p-q)/q},
\end{equation}
where the supremum is taken over all sequences of sets as above.
The following result is a metric space version of a theorem
on page 120 of \cite{mazja85}.

\begin{theorem}\label{theorem:qp}
  Suppose that $1\leq q<p<\infty$.
  \begin{itemize}
  \item[(i)] If $\gamma<\infty$, then
    \begin{equation}\label{eqn:sp}
      \left(\int_\Omega |u|^q\dd\nu\right)^{1/q}
      \leq c_S\left(\int_\Omega g_u^p\dd\mu\right)^{1/p},
    \end{equation}
    where $c_S$ is independent of $u\in N^{1,p}_0(\Omega)\cap C(\Omega)$.
  \item[(ii)] If~\eqref{eqn:sp} holds for every $u\in N^{1,p}_0(\Omega)\cap C(\Omega)$ 
    and if the constant $c_S$ is independent of $u$, then $\gamma<\infty$.
  \end{itemize}
\end{theorem}

\begin{proof}
  (i)
  We have 
  \begin{equation*}
    \begin{split}
      \int_\Omega|u|^q\dd\nu=
      &\sum_{j=-\infty}^\infty 
      q\int_{2^{j}}^{2^{j+1}}\lambda^{q-1}\nu(E_\lambda)\dd\lambda\\
      \leq&\, q\sum_{j=-\infty}^\infty 
      2^{(j+1)(q-1)}2^j\nu(E_{2^j})\\
      =&\, q 2^{q-1}\sum_{j=-\infty}^\infty 2^{jq}\nu(E_{2^j}).
    \end{split}
  \end{equation*}
  By the H\"older inequality,
  \begin{equation*}
    \begin{split}
      \sum_{j=-\infty}^\infty& 2^{jq}\nu(E_{2^j})\\
      =& \sum_{j=-\infty}^\infty
      \left(
        \frac{\nu(E_{2^j})^{p/q}}
        {\Ca_p( E_{2^j},E_{2^{j-1}})}\right)^{q/p}
      \left(
        2^{jp}\Ca_p( E_{2^j},E_{2^{j-1}})
      \right)^{q/p}\\
      \leq&  
      \left( \sum_{j=-\infty}^\infty  
        \left( \frac{\nu(E_{2^j})^{p/q}}
          {\Ca_p( E_{2^j},E_{2^{j-1}})}
        \right)^{q/(p-q)}  
      \right)^{(p-q)/p}
      \\   
      &\times\left( \sum_{j=-\infty}^\infty 
        2^{jp}\Ca_p(E_{2^j},E_{2^{j-1}})  
      \right)^{q/p}.
    \end{split}
  \end{equation*}   
  Let
  \begin{equation*}
      u_j=\begin{cases}
        1, &\textrm{if}\quad|u|> 2^j,\\
        \dfrac{|u|-2^{j-1}}{2^{j-1}},
        &\textrm{if}\quad 2^{j-1}<|u|\leq 2^j,\\
        0,&\textrm{if}\quad|u|\leq 2^{j-1}.
      \end{cases}
    \end{equation*}
    Then
    \begin{equation*}
      \begin{split}
        \Ca_p(E_{2^j},E_{2^{j-1}})\leq 
        \int_\Omega g_{u_j}^p\dd\mu
        \leq  2^{-(j-1)p}
        \int_{E_{2^{j-1}}\setminus E_{2^j}}g_u^p\dd\mu
      \end{split}
    \end{equation*}
    It follows that
    \begin{equation*}
      \begin{split}
      \sum_{j=-\infty}^\infty 
      2^{jp}\Ca_p(E_{2^j},E_{2^{j-1}})
      &\leq 2^p \sum_{j=-\infty}^\infty 
      \int_{E_{2^{j-1}}\setminus E_{2^j}}g_u^p\dd\mu\\
      &= 2^p\int_\Omega g_u^p\dd\mu
     \end{split}
    \end{equation*}
    and consequently
    \[
    \int_\Omega |u|^q\dd\nu\leq c\left(\int_\Omega g_u^p\dd\mu\right)^{q/p}.
    \]
    
(ii)
Let $E_j$ be as as in the statement of the theorem, and define
    \[ 
      \lambda_j=\displaystyle{\sum_{i=j}^N
      \left(\frac{\nu(E_i)}{\Ca_p( E_i,E_{i+1})}
      \right)^{1/(p-q)}},\quad j=-N,-N+1,\ldots, N,
     \]
and $\lambda_{N+1}=0$. 
    Let $u_j\in\mathcal A(E_j,E_{j+1})$ be continuous, and define
    \begin{equation*}
      u=
      \begin{cases}
        (\lambda_j-\lambda_{j+1})u_j+\lambda_{j+1} 
        &\textrm{in}\quad E_{j+1}\setminus E_j,\\
        \lambda_{-N}&\textrm{in}\quad E_{-N},\\
        0&\textrm{in}\quad\Omega\setminus E_{N+1}.
      \end{cases}
    \end{equation*}
    Then $u\in N^{1,p}_0(\Omega)\cap C(\Omega)$. By the Cavalieri principle
    \begin{equation*}
      \begin{split}
        \int_\Omega|u|^q\dd\nu =
        &q\int_0^\infty\lambda^{q-1}\nu(E_\lambda)\dd\lambda
        =\sum_{j=-N}^Nq\int_{\lambda_{j+1}}^{\lambda_j}
        \lambda^{q-1}\nu(E_\lambda)\dd\lambda\\
        \geq&\sum_{j=-N}^N\nu(E_j)(\lambda_j^q-\lambda_{j+1}^q).
      \end{split}
    \end{equation*}
    From this we conclude that
    \begin{equation*}
      \begin{split}
        \bigg(\sum_{j=-N}^N&\nu(E_j)(\lambda_j-\lambda_{j+1})^q 
        \bigg)^{p/q}
        \leq\left( \sum_{j=-N}^N \nu(E_j) (\lambda_j^q-\lambda_{j+1}^q) 
        \right)^{p/q}\\
        &\leq\left( \int_\Omega |u|^q \dd\nu\right)^{p/q} 
        \leq c_S\int_\Omega g_u^p \dd\nu
        =c_S\sum_{j=-N}^N \int_{E_{j+1}\setminus E_j} g_u^p\dd\mu\\
        &\leq c_S\sum_{j=-N}^N (\lambda_j-\lambda_{j+1})^p 
        \int_{E_{j+1}\setminus E_j} g_{u_j}^p\dd\mu.
      \end{split}
    \end{equation*}
    Taking the infimum on the right-hand side, we arrive at
    \begin{equation*}
      \left(\sum_{j=-N}^N\nu(E_j)(\lambda_j-\lambda_{j+1})^q\right)^{p/q}
      \leq c_S\sum_{j=-N}^N(\lambda_j-\lambda_{j+1})^p
      \Ca_p(E_j,E_{j+1}).
    \end{equation*}
    Since
    \[
    \lambda_j-\lambda_{j+1}
    =\left(\frac{\nu(E_j)}{\Ca_p(E_j,E_{j+1})}
    \right)^{1/(p-q)},
    \]
    we obtain
    \[
    \begin{split}
      &\left[        
        \sum_{j=-N}^N         
        \left(           \frac{\nu(E_j)^{p/q}}{\Ca_p(E_j,E_{j+1})}    
        \right)^{q/(p-q)} 
      \right]^{p/q}\\
      &\qquad\qquad\leq
      c_S\sum_{j=-N}^N 
      \left( 
        \frac{\nu(E_j)^{p/q}}{\Ca_p(E_j,E_{j+1})} 
      \right)^{q/(p-q)},
    \end{split}
    \]
    and the claim follows.
\end{proof}

Next we present an integral version of Theorem~\ref{theorem:qp}. See
also page 30 in~\cite{mazja03} for the Euclidean case.

\begin{theorem}
Let $ q< p$ and $\mu(\Omega)<\infty$ and
\[
\lambda_p(s)=\inf \{\Ca_p(G)\,:\,G\subset\Omega\textrm{ and }\nu(G)\geq s\}.
\]
Then
\begin{equation}\label{eqn:integraali}
  \int_0^{\nu(\Omega)}\left(\frac{t^{p/q}}{\lambda_p(t)}\right)^{q/(p-q)}\frac{\dd t}{t}\leq c_I<\infty
\end{equation}
if and only if the Sobolev inequality~\eqref{eqn:sp} holds for every
$u\in N_0^{1,p}(\Omega)\cap C(\Omega)$ with a constant $c_S$ that is
independent of $u$.
\end{theorem}

\begin{proof}
  First, assume that~\eqref{eqn:integraali} holds. Let $s=p/q$ and
  $s'=p/(p-q)$ be the H\"older conjugate of $s$.  Then

\begin{equation*}
\begin{split}
  \int_\Omega u^q&\dd\nu \leq  2^q\sum_{j= -\infty}^\infty 2^{jq}\nu(\{2^j<u<2^{j+1}\})\\
  &\leq  2^q\left(\sum_{j=-\infty}^\infty 2^{jp}\Ca_p(\{u>2^j\})\right)^{q/p}\\
  &\qquad\times \left( \sum_{j=-\infty}^\infty \frac{ (\nu( \{ u >
      2^j\} )-\nu( \{u > 2^{j+1} \} ))^{s'}}{\Ca_p( \{u > 2^j\}
      )^{s'/s} }
  \right)^{1/s'}\\
  &\leq c\left(\int_\Omega g_u^p\dd\mu\right)^{q/p}
  \left(\int_0^{\nu(\Omega)} \frac{t^{s'-1}}{\lambda_p(t)^{s'/s}}\dd t\right)^{1/s'}\\
&  \leq c_I^{1/s'}c\left(\int_\Omega g_u^p\dd\mu\right)^{q/p}.
\end{split}
\end{equation*}

Here we used the H\"older inequality, monotonicity of $\lambda_p(t)$ and the fact that
\[
\Ca_p(\{u>2^j\})\leq 2^{-jp+p}\int_{\{2^{j-1}<u<2^j\}}g_u^p\dd\mu.
\]
Assume now that~\eqref{eqn:sp} holds. For every $j\in \Z$, let $u_j\in
N_0^{1,p}(\Omega)\cap C(\Omega)$ be a function such that $0\leq
u_j\leq 1$, $\nu(\{u_j=1\})\geq 2^j$ and
\[
\int_\Omega g_{u_j}^p\dd\mu\leq \lambda_p(2^j)+\varepsilon_j,
\]
with $
0\leq\varepsilon_j\leq\lambda_p(2^j)$.
Let
\[
u=\sup_j \beta_j u_j,
\]
where
\[
\beta_j=\left(\frac{2^j}{\lambda_p(2^j)}\right)^{1/(p-q)}.
\]
Now
\begin{equation}\label{eqn:fq}
\int_\Omega u^q \dd\nu\geq \frac12\sum_{j=-\infty}^\infty \beta_j^q 2^j
\end{equation}
and
\begin{equation}\label{eqn:gp}
\int_\Omega g_u^p\dd\mu\leq \left(\sum_j \beta_j^p(\lambda_p(2^j)+\varepsilon_j)\right)\leq 2\sum_j \beta_j^p\lambda_p(2^j).
\end{equation}
As
\[
\beta_j^q 2^j=\beta_j^p\lambda_p(2^j),
\]
it follows by~\eqref{eqn:sp},~\eqref{eqn:fq} and~\eqref{eqn:gp}  that
\[
\sum_j \frac{(2^j)^{p/(p-q)}}{\lambda_p(2^j)^{q/(p-q)}}=\sum_j \beta_j^q 2^j\leq c,
\]
and
\[
\int_0^{\nu(\Omega)}\left(\frac{t^{p/q}}{\lambda_p(t)}\right)^{q/(p-q)}\frac{\dd t}{t}\leq c.
\]
by monotonicity of $\lambda_p$.\end{proof}

\section{Functions with general boundary values}
In this section, we obtain necessary and sufficient conditions for 
Sobolev--Poincar\'e inequalities of type
\[
\inf_{a\in\R}
\left(
\int_\Omega |u-a|^q\dd\nu
\right)^{1/q}
\leq c\left(\int_\Omega g_u^p\dd\mu\right)^{1/p}
\]
where $u\in N^{1,p}(\Omega)\cap C(\Omega)$ and $1\leq p\leq q<\infty$.
To this end, we shall need the concept of conductivity, see Chapter 4
in \cite{mazja85}.

Let $\Omega\subset X$ be a bounded open set. Let $F$ be a closed
subset of $\Omega$ and let $G$ be an open subset of $\Omega$ such that
$F\subset G$. The open set $C=G\setminus F$ is called a
\emph{conductor} and
\[
\mathcal B(F,G,\Omega)=\{u\in N^{1,p}(\Omega)\cap C(\Omega)\,:\, u\geq
1\textrm{ in } F\textrm{ and } u\leq 0\textrm{ in }\Omega\setminus G\}
\]
is the set of admissible functions. The number
\[
\con_p(F,G,\Omega)=\inf_{u\in \mathcal B(F,G,\Omega)}\int_\Omega
g_u^p\dd\mu
\]
is called the \emph{$p$--conductivity} of $C$.

The next result can be proved in the same way as Lemma~\ref{lemma:capaint}.

\begin{lemma}\label{lemma:conint}
  Let $G\subset \Omega$ be open and $1\leq p<\infty $. Suppose that 
  $u\in N^{1,p}(\Omega)\cap C(\Omega)$ such that $u=0$ in $\Omega\setminus G$. Then
  \[
  \int_0^\infty\lambda^{p-1}\con_p(E_\lambda,G,\Omega)
  \dd\lambda\leq 2^{2p-1}\int_\Omega g_u^p\dd\mu.
  \]
\end{lemma}

The following result is a metric space version of a theorem
on page 210 of \cite{mazja85}.

\begin{theorem}\label{theorem:theorem2}
  Suppose that $1\leq p\leq q<\infty$ and let $G\subset \Omega$ be open.
  \begin{itemize}
  \item[(i)] If there is a constant $\gamma$ such that
    \begin{equation}\label{eqn:nu_con}
      \nu(F)^{p/q}\leq \gamma \con_p(F,G,\Omega)
    \end{equation}
    for every $F\subset G$, then
    \begin{equation}\label{eqn:nu_g}
      \left(\int_\Omega|u|^q\dd\nu\right)^{1/q}\leq
      c\left(\int_\Omega g_u^p\dd\mu\right)^{1/p}
    \end{equation}
    for every $u\in N^{1,p}(\Omega)\cap C(\Omega)$ such that $u=0$ in $\Omega\setminus G$. 
    Here the constant $c$ is independent of $u$.
  \item[(ii)] If~\eqref{eqn:nu_g} holds for every $u\in N^{1,p}(\Omega)\cap C(\Omega)$ 
    such that $u=0$ in $\Omega\setminus G$, then~\eqref{eqn:nu_con} holds 
    for every $F\subset G$ with $\gamma= c$.
  \end{itemize}
\end{theorem}

\begin{proof}

(i)
    We conclude by Lemma~\ref{lemma:qpey},~\eqref{eqn:nu_con} and 
    Lemma~\ref{lemma:conint} that
    \begin{equation*}
      \begin{split}
        \left(\int_\Omega |u|^q\dd\nu\right)^{1/q}\leq&
        \left(p\int_0^\infty\lambda^{p-1}\nu(E_\lambda)^{p/q}
          \dd\lambda\right)^{1/p}\\
        \leq&\left(p\gamma\int_0^\infty\lambda^{p-1}\con_p(E_\lambda, G,\Omega)
          \dd\lambda\right)^{1/p}\\
        \leq&\left(p\gamma 2^{2p-1}\int_\Omega g_u^p\dd\mu\right)^{1/p}.
      \end{split}
    \end{equation*}
    
(ii)
    If $u\in \mathcal B(F,G,\Omega)$, then
    inequality~\eqref{eqn:nu_g} implies that
    \begin{equation*}
      \begin{split}
        \nu(F)^{1/q}\leq&\left(\int_\Omega|u|^q\dd\nu\right)^{1/q}
        \leq c\left(\int_\Omega g_u^p\dd\mu\right)^{1/p}.
      \end{split}
    \end{equation*}
    The claim follows by taking the infimum on the right-hand side.
\end{proof}

\begin{theorem}
  Suppose that $1\leq p\leq q<\infty$ and that $\Omega\subset X$ 
  is a bounded open set.
  \begin{itemize}
  \item[(i)] If there is a constant $\gamma$ such that
    \begin{equation}\label{eqn:r_nu_con}
      \nu(F)^{p/q}\leq \gamma\con_p(F,G,\Omega)
    \end{equation}
    for every conductor $G\setminus F$ with $\nu(G)\leq
    \nu(\Omega)/2$, then
    \begin{equation}\label{eqn:r_nu_g}
      \inf_{a\in\R}\left(
        \int_\Omega |u-a|^q\dd\nu
      \right)^{1/q}
      \leq c\left(\int_\Omega g_u^p\dd\mu\right)^{1/p}
    \end{equation}
for every $u\in N^{1,p}(\Omega)\cap C(\Omega)$.
\item[(ii)] If\eqref{eqn:r_nu_g} holds for every $u\in
  N^{1,p}(\Omega)\cap C(\Omega)$, then~\eqref{eqn:r_nu_con} holds for
  every conductor $G\setminus F$ with $\nu(G)\leq \nu(\Omega)/2$.
  \end{itemize}
\end{theorem}

\begin{proof}

(i)
Let $\alpha\in\R$ be such that
    \[
     \nu(\{x\in\Omega\,:\, u(x)\geq \alpha\})\geq \frac12\nu(\Omega) 
      \]
and
\[
\nu(\{x\in\Omega\,:\, u(x)> \alpha\})\leq \frac12\nu(\Omega).
\]
Now $(u-\alpha)_+\in N^{1,p}(\Omega)\cap C(\Omega)$ and $u=0$ 
in $\Omega\setminus G$, where
\[
  G=\{x\in\Omega\,:\, u(x)>\alpha\}.
\]
Clearly $\nu(G)\leq \frac12 \nu(\Omega)$ and by~\eqref{eqn:r_nu_con}, 
\[\nu(F)^{p/q}\leq\gamma\con_p(F,G,\Omega)
\] 
for every $F\subset G$.
By Theorem~\ref{theorem:theorem2}, we have
\[
\left(\int_\Omega (u-\alpha)_+^q\dd\nu\right)^{1/q}\leq c 
\left(\int_{\{x\in\Omega\,:\, u(x)>\alpha\}} g_u^p\dd\mu\right)^{1/p}.
\]
Similarly,
\[
\left(\int_\Omega (\alpha-u)_+^q\dd\nu\right)^{1/q}
\leq c\left(
  \int_{\{x\in\Omega\,:\, u(x)<\alpha\}} g_u^p\dd\mu\right)^{1/p}.
\]
A combination of these estimates implies that
\begin{equation*}
  \begin{split}
    \bigg( \int_\Omega&|u-\alpha|^q\dd\mu \bigg)^{1/q}\\ 
    &\leq \left( \int_\Omega (\alpha-u)_+^q \dd\mu\right)^{1/q} 
    + \left( \int_\Omega (u-\alpha)_+^q \dd\mu\right)^{1/q}\\
    &\leq c\left(\int_{\{x\in\Omega\,:\, u(x)>\alpha\}}g_u^p\dd\mu
    \right)^{1/p} 
    +c\left(\int_{\{x\in\Omega\,:\, u(x)<\alpha\}}g_u^p\dd\mu
    \right)^{1/p}\\
    &\leq c\left(\int_\Omega g_u^p\dd\mu\right)^{1/p}.
  \end{split}
\end{equation*}

(ii) 
Let $G\setminus F$ be a conductor with $\nu(G)\leq\frac12\nu(\Omega)$ 
and suppose that $u\in N^{1,p}(\Omega)\cap C(\Omega)$ such that $u=0$ 
in $\Omega\setminus G$ and $u=1$ on $F$. Since
    \begin{equation*}
      \begin{split}
        \bigg(\int_\Omega &|u-u_\Omega|^q\dd\nu
        \bigg)^{1/q}\\
        &\leq \left(\int_\Omega |u-a|^q\dd\nu \right)^{1/q} 
        + |a-u_\Omega|\nu(\Omega)^{1/q}\\
        &\leq \left(\int_\Omega |u-a|^q\dd\nu \right)^{1/q} 
        +\nu(\Omega)^{1/q}\kint_\Omega |u-a|\dd\nu\\
        &\leq \left(\int_\Omega |u-a|^q\dd\nu \right)^{1/q} 
        +\nu(\Omega)^{1/q}\left(\kint_\Omega|u-a|^q\dd\nu\right)^{1/q}\\
        &= 2\left(\int_\Omega |u-a|^q\dd\nu\right)^{1/q},
      \end{split}
    \end{equation*}
    we have
    \begin{equation*}
     \begin{split}
      \inf_{a\in\R} \left(\int_\Omega |u-a|^q\dd\nu\right)^{1/q}
      &\leq \left(\int_\Omega |u-u_\Omega|^q\dd\nu\right)^{1/q}\\
      &\leq 2\inf_{a\in\R} \left(\int_\Omega |u-a|^q\dd\nu\right)^{1/q}.
     \end{split}
    \end{equation*}
    Now
    \begin{equation*}
      \begin{split}
        c\left(\int_\Omega g_u^p\dd\mu\right)^{q/p}
         &\geq\int_\Omega|u-u_\Omega|^q\dd\nu\\
         &=\int_G|u-u_\Omega|^q\dd\nu + |u_\Omega|^q\nu(\Omega\setminus G).
      \end{split}
    \end{equation*}
    This and the fact that 
    $\nu(\Omega\setminus G)\geq\frac{1}{2}\nu(\Omega)$ imply that
    \[
    |u_\Omega|^q\nu(\Omega)
    \leq c\left(\int_\Omega g_u^p\dd\mu\right)^{q/p}.
    \]
    Since
    \begin{equation*}
      \begin{split}
      \left(\int_\Omega |u|^q\dd\nu\right)^{1/q}
      &\leq |u_\Omega|\nu(\Omega)^{1/q}
      + \left(\int_\Omega |u-u_\Omega|^q\dd\nu\right)^{1/q}\\
      &\leq c\left(\int_\Omega g_u^q\dd\mu\right)^{1/p},
      \end{split}
    \end{equation*}
    we have
    \begin{equation*}
      \nu(F)^{1/q}
      \leq\left( \int_\Omega |u|^q\dd\nu \right)^{1/q}
      \leq c\left( \int_\Omega g_u^p\dd\mu\right)^{1/p}
    \end{equation*}
    and by taking an infimum over all functions $u$, we have
    \[
    \nu(F)^{p/q}\leq c\con_p(F,G,\Omega).
    \] 
\end{proof}

\def\cprime{$'$} \def\ocirc#1{\ifmmode\setbox0=\hbox{$#1$}\dimen0=\ht0
  \advance\dimen0 by1pt\rlap{\hbox to\wd0{\hss\raise\dimen0
  \hbox{\hskip.2em$\scriptscriptstyle\circ$}\hss}}#1\else {\accent"17 #1}\fi}

\vspace{0.5cm}
\noindent
\small{\textsc{J.K.},}
\small{\textsc{Institute of Mathematics},}
\small{\textsc{P.O. Box 1100},}
\small{\textsc{FI-02015 Helsinki University of Technology},}
\small{\textsc{Finland}}\\
\footnotesize{\texttt{juha.kinnunen@tkk.fi}}

\vspace{0.3cm}
\noindent
\small{\textsc{R.K.},}
\small{\textsc{Institute of Mathematics},}
\small{\textsc{P.O. Box 1100},}
\small{\textsc{FI-02015 Helsinki University of Technology},}
\small{\textsc{Finland}}\\
\footnotesize{\texttt{rkorte@math.tkk.fi}}

\end{document}